# On the role of Riesz potentials in Poisson's equation and Sobolev embeddings

Rahul Garg
rgarg@tx.technion.ac.il

Daniel Spector
dspector@tx.technion.ac.il

Department of Mathematics
Technion - Israel Institute of Technology

**Abstract**

In this paper, we study the mapping properties of the classical Riesz potentials acting on $L^p$-spaces. In the supercritical exponent, we obtain new "almost" Lipschitz continuity estimates for these and related potentials (including, for instance, the logarithmic potential). Applications of these continuity estimates include the deduction of new regularity estimates for distributional solutions to Poisson's equation, as well as a proof of the supercritical Sobolev embedding theorem first shown by Brezis and Wainger in 1980.

# Contents

# 1 Introduction and Statement of Results

For each $0 < \alpha < N$, one defines the Riesz potential

$$I_\alpha(x) := \frac{\gamma(N, \alpha)}{|x|^{N-\alpha}}.$$

Then the operator $I_\alpha$ can be defined through its action on suitable subspaces of the space of Lebesgue measurable functions by convolution in the sense that

$$I_\alpha f(x) := (I_\alpha * f)(x) = \gamma(N, \alpha) \int_{\mathbb{R}^N} \frac{f(y)}{|x - y|^{N-\alpha}} \, dy, \qquad (1.1)$$

whenever it is well-defined. Here, the precise definition of the constant is

$$\gamma(N,\alpha) := \frac{\Gamma(\frac{N-\alpha}{2})}{\pi^{\frac{N}{2}} 2^\alpha \Gamma(\frac{\alpha}{2})},$$

from which one obtains two key properties of the potentials: the semigroup property and its connection with the Laplacian (see, for example, Stein [20][Chapter 5, p.118]).

In the next section, we review the literature concerning the mapping properties of the Riesz potentials on $L^p(\mathbb{R}^N)$ for $0 < \alpha < N$ and $1 < p \leq \frac{N}{\alpha-1}$. As we will see, though the literature is quite extensive, notably absent is a complete treatment of the supercritical case $p = \frac{N}{\alpha-1}$, as well as the consideration of $\alpha \geq N$. Although these might at first seem like purely academic questions, a closer examination reveals connections between these considerations and several interesting applications. Indeed, in the supercritical case $\alpha = 1 + \frac{N}{p}$, $\alpha \in [1, N)$, we show that one has "almost" Lipschitz continuity estimates for the Riesz potentials, while we introduce some Riesz-type potentials defined for $\alpha \in [N, N+1)$ and demonstrate that they enjoy analogous regularity properties to the Riesz potentials. Interesting in its own right are the new estimates we show for the logarithmic potential (which corresponds to the case $\alpha = N$), though more generally we demonstrate some important consequences of our results through the connection of potentials, partial differential equations, and Sobolev embeddings. This includes the conclusion of new regularity results for Poisson's equation in any number of dimensions, as well as a potential approach to the Sobolev embedding theorem in the supercritical exponent. In particular, this shows how Riesz and Riesz-type potentials enable one to give a unified treatment of Sobolev embeddings for $p > 1$. We now proceed to review the literature concerning the Riesz potentials, after which we will state our results in Section 1.2.

## 1.1 Historical Background

The study of the mapping properties of $I_\alpha$ in several dimensions was initiated by S. Sobolev in his 1938 paper [19], where he proved the following theorem concerning the $L^{p^*}$-summability of $I_\alpha f$ for functions $f \in L^p(\mathbb{R}^N)$ with $1 < p < \frac{N}{\alpha}$.[1]

**Theorem 1.1** *Let $0 < \alpha < N$ and $1 < p < \frac{N}{\alpha}$. Then there exists $C = C(p, \alpha, N)$ such that*

$$\|I_\alpha f\|_{L^{p^*}(\mathbb{R}^N)} \leq C\|f\|_{L^p(\mathbb{R}^N)}$$

*for all $f \in L^p(\mathbb{R}^N)$, where $p^* := \frac{Np}{N-\alpha p}$.*

When $p = \frac{N}{\alpha}$, one cannot expect an embedding theorem of the above type, and the appropriate replacement is exponential integrability of the potential. When $p = N = 2$, the following theorem is related to Lemma 1 in Pohožaev's 1965 paper [17], while more generally the potential arguments are due to Trudinger's 1967 paper [23] for $\alpha = 1$. For general $\alpha$, the result was obtained by Strichartz [22] and also Hedberg [8].[2] Here we record the statement of the theorem as can be found in the book of Mizuta [11][Section 4.2, Theorem 2.3, p. 156].

---
[1]The reference of Stein [20][Chapter V, p.165] gives some details concerning Sobolev's contribution in the context of the one dimensional result of Hardy and Littlewood.

**Theorem 1.2** *Let $0 < \alpha < N$ and $p = \frac{N}{\alpha}$. Then for every bounded, open subset $G \subset \mathbb{R}^N$, there exist $A_1, A_2 > 0$ such that*

$$\frac{1}{|G|} \int_G exp\left(\left[\frac{|I_\alpha f(x)|}{A_1 \|f\|_{L^p(\mathbb{R}^N)}}\right]^{p'}\right) dx \leq A_2$$

*for all $f \in L^p(G)$.*

An early reference which obtains a partial result of the preceding theorem is the 1961 paper of Yudovič [24], which also records a version of the following theorem concerning the case $\frac{N}{\alpha} < p < \frac{N}{\alpha-1}$, assuming $f$ has support in a bounded domain. We again give a statement of this theorem from the book of Mizuta (see [11][Section 4.2, Theorem 2.2, p. 155]).

**Theorem 1.3** *Let $1 < p < +\infty$ and $0 < \alpha < N$ be such that $0 < \alpha - \frac{N}{p} < 1$. Then there exists $C = C(p, \alpha, N)$ such that*

$$|I_\alpha f(x) - I_\alpha f(z)| \leq C|x-z|^{\alpha - \frac{N}{p}} \|f\|_{L^p(\mathbb{R}^N)}$$

*for all $f \in L^p(\mathbb{R}^N)$ such that $I_\alpha f$ is well-defined.*

**Remark 1.4** *Let us remark that for this range of $p$ and $\alpha$, it is not the case that $I_\alpha f$ is well-defined for every $f \in L^p(\mathbb{R}^N)$. By considering functions with support in a bounded domain, one has $f \in L^q(\mathbb{R}^N)$ for every $1 \leq q \leq \frac{N}{\alpha}$, so that by the preceding results one has $I_\alpha f$ is well-defined. Working on all of the space, we instead prefer to introduce a modified Riesz potential*

$$\tilde{I}_\alpha f(x) := \gamma(N, \alpha) \int_{\mathbb{R}^N} f(y) \left[\frac{1}{|x-y|^{N-\alpha}} - \frac{1}{|y|^{N-\alpha}}\right] dy,$$

*since $\tilde{I}_\alpha f$ is well-defined for all $f \in L^p(\mathbb{R}^N)$ with $0 < \alpha - \frac{N}{p} < 1$. Moreover,*

$$I_\alpha f(x) - I_\alpha f(z) = \tilde{I}_\alpha f(x) - \tilde{I}_\alpha f(z),$$

*whenever $I_\alpha f$ is well-defined.*

Finally, when $p = \frac{N}{\alpha-1}$, Yudovič [24] had announced the following theorem for potentials acting on functions with compact support.

**Theorem 1.5** *Suppose $\Omega \subset \mathbb{R}^N$ is open and bounded. Let $1 < p < +\infty$ and $0 < \alpha < N$ be such that $p = \frac{N}{\alpha-1}$. Then there exists $C = C(p, \alpha, N, \Omega)$ such that*

$$|I_\alpha f(x) - I_\alpha f(z)| \leq C|x-z|\left(|\ln|x-z|| + 1\right)^{\frac{1}{p'}} \|f\|_{L^p(\Omega)}$$

*for all $f \in L^p(\mathbb{R}^N)$ such that $supp\ f \subset \Omega$.*

---

[2]Strichartz's paper actually concerns the Bessel potential, and so for our purposes Hedberg is the more appropriate reference.

Aside from the question of the removal of the assumption of compact support in Theorem 1.5, this gives a complete picture on the mapping properties of the Riesz potentials for $1 < p \leq \frac{N}{\alpha-1}$. When $\frac{N}{\alpha-1} < p$, observing this requires $\alpha > 1$, we find that $I_\alpha f$ has a distributional derivative which is essentially $I_{\alpha-1}f$. Thus, one can apply the preceding theorems to the derivatives of the potential and obtain a family of results for $p$ in this range (for a fixed $\alpha$). There are also a considerable number of results for Riesz potentials beyond the basic framework of the $L^p$-spaces, for example, the consideration of potentials with variable exponent [13], potentials mapping on $L^p$-spaces with variable exponent [3], potentials acting on functions where the underlying space is assumed to be metric [4], potentials acting on Morrey spaces of variable exponent [14], potentials acting on functions in general Orlicz spaces [15] or Musielak-Orlicz-Morrey spaces [16]. This list is by no means exhaustive, though gives an idea of some of the possible variations one can consider and obtain results analogous to those we have recorded here.

## 1.2 Statement of Results

Analyzing the literature, several natural questions present themselves for study. Firstly, one wonders whether it is possible to remove the hypothesis of compact support from Theorem 1.5. As we can see in the following theorem, one can relax this hypothesis by asking the function decay sufficiently to have some global summability in another $L^q$-space.

**Theorem A** *Let $1 < p < +\infty$ and $\alpha = 1 + \frac{N}{p} < N$. Then for any $1 \leq q < p$, there exists $C = C(q, \alpha, N)$ such that*

$$|\tilde{I}_\alpha f(x) - \tilde{I}_\alpha f(z)| \leq C|x-z|(|\ln|x-z||+1)^{\frac{1}{p'}}\left(\|f\|_{L^q(\mathbb{R}^N)} + \|f\|_{L^p(\mathbb{R}^N)}\right)$$

*for all $f \in L^p(\mathbb{R}^N) \cap L^q(\mathbb{R}^N)$.*

As alluded to earlier in the introduction, another interesting question is that of the restriction $0 < \alpha < N$ in the definition (1.1). Whereas in the case $\alpha \leq 0$, the potentials $I_\alpha(x) = \gamma(N,\alpha)/|x|^{N-\alpha}$ cease to be integrable in a neighborhood of the origin, when $\alpha > N$, $I_\alpha$ does not have a local singularity. For example, the fundamental solution of Laplace's equation in the plane is known to be the logarithmic potential, and more generally one can in higher dimensions define the logarithmic potential[3]

$$I_N f(x) := \frac{\Gamma(\frac{N}{2})}{2\pi^{\frac{N}{2}}} \int_{\mathbb{R}^N} \log \frac{1}{|x-y|} f(y)\,dy. \qquad (1.2)$$

A comparison of $I_N$ and $I_\alpha$ reveals that the known continuity estimates for the logarithmic potential are not sharp when compared with those known for the Riesz potentials, even in **the same regime of criticality**. For example, the estimates obtained in Theorem 1.3 and Theorem A are not stable under the limit $\alpha \to N$, and therefore yield no useful information regarding the logarithmic potential. Moreover, as we will see in Section 2, the naive approach to mimic the techniques from Riesz potential estimates also do not obtain the sharp exponent. Nonetheless, we are able to prove the analogous theorem concerning the continuity properties of the modified logarithmic potential (which can be compared, for instance, with [10][Chapter 10, p. 260]).

---

[3]The logarithmic potential can also be obtained as a limiting process as $\alpha \to N$ for certain classes of functions, see for example, [9][Chapter I, p. 50].

**Theorem B** *Let $1 < p \leq \frac{N}{N-1}$.*

*i) If $1 < p < \frac{N}{N-1}$, then there exists $C = C(p, N)$ such that*

$$|\tilde{I}_N f(x) - \tilde{I}_N f(z)| \leq C|x - z|^{N - \frac{N}{p}} \|f\|_{L^p(\mathbb{R}^N)}$$

*for all $f \in L^p(\mathbb{R}^N)$.*

*ii) If $p = \frac{N}{N-1}$ and $1 \leq q < p$, then there exists $C = C(q, N)$ such that*

$$|\tilde{I}_N f(x) - \tilde{I}_N f(z)| \leq C|x - z| \left(|\ln|x - z|| + 1\right)^{\frac{1}{N}} \left(\|f\|_{L^p(\mathbb{R}^N)} + \|f\|_{L^q(\mathbb{R}^N)}\right)$$

*for all $f \in L^p(\mathbb{R}^N) \cap L^q(\mathbb{R}^N)$.*

Here, in analogy to the modified Riesz potentials, we have considered the modified logarithmic potential as follows.

**Definition 1.6** *We define the modified Logarithmic potential*

$$\tilde{I}_N f(x) := \frac{\Gamma(\frac{N}{2})}{2\pi^{\frac{N}{2}}} \int_{\mathbb{R}^N} \left[\log \frac{1}{|x - y|} - \log \frac{1}{|y|}\right] f(y) \, dy. \tag{1.3}$$

Since for $N \geq 2$ there is an intimate relationship between certain distributional solutions to Poisson's equation and the potentials $I_2$, an immediate application of the continuity results asserted in Theorem A and Theorem B is to conclude sharp and uniform regularity estimates for these solutions (such a results has been announced in [6]). We now make this statement more precise.

We suppose that $f \in L^p(\mathbb{R}^N)$ for some $\frac{N}{2} < p < N$, and we say that $u \in L^1_{loc}(\mathbb{R}^N)$ is a distributional solution to Poisson's equation if

$$-\int_{\mathbb{R}^N} u \Delta \varphi = \int_{\mathbb{R}^N} f \varphi, \tag{1.4}$$

for all $\varphi \in C_c^\infty(\mathbb{R}^N)$. Now, we know that there exists $u \in L^1_{loc}(\mathbb{R}^N)$ which satisfies (1.4) (see, for example, [10][Chapter 6, Theorem 6.21, p.157]), and that if $supp\, f$ is compact, this solution $u$ is given by $u = I_2 f$. Moreover, any two solutions differ only by a harmonic function, and therefore if we restrict ourselves to solutions which satisfy the growth condition

$$\frac{u(x)}{|x|} \to 0 \text{ as } |x| \to \infty, \tag{1.5}$$

then, up to a constant, such a $u \in L^1_{loc}(\mathbb{R}^N)$ is unique. Then we have the following theorem on the regularity of this solution.

**Theorem C** *Suppose $N \geq 2$ and $\frac{N}{2} < p \leq N$. If $\frac{N}{2} < p < N$, assume $f \in L^p(\mathbb{R}^N)$, while if $p = N$ assume $f \in L^N(\mathbb{R}^N) \cap L^q(\mathbb{R}^N)$ for some $1 < q < N$. Then $u = \tilde{I}_2 f$ satisfies (1.4) and has the following regularity estimates:*

i) If $\frac{N}{2} < p < N$, then
$$|u(x) - u(z)| \leq C|x-z|^{2-\frac{N}{p}} \|f\|_{L^p(\mathbb{R}^N)}.$$

ii) If $p = N$, then
$$|u(x) - u(z)| \leq C|x-z| (|\ln|x-z|| + 1)^{\frac{1}{N'}} \left( \|f\|_{L^N(\mathbb{R}^N)} + \|f\|_{L^q(\mathbb{R}^N)} \right).$$

Our approach to Theorem B, and therefore the mechanism behind the better regularity results deduced in Theorem C, is a new representation for the logarithmic potential. In fact, this representation is one way to extend the Riesz potentials past the criticality $\alpha = N$ to all $\alpha \in [N, N+1)$.

**Definition 1.7** Let $\beta \in [0,1)$ and define
$$\tilde{T}_j^\beta f(x) := \frac{\gamma(N, N-1+\beta)}{N-1+\beta} \int_{\mathbb{R}^N} \left[ \frac{y_j - x_j}{|y-x|^{1-\beta}} - \frac{y_j}{|y|^{1-\beta}} \right] f(y) \, dy. \quad (1.6)$$

We then have the following theorem concerning the mapping properties of $\tilde{T}_j^\beta$.

**Theorem D** Let $\beta \in [0,1)$.

i) If $1 < p < \frac{N}{N-1+\beta}$, then there exists $C = C(p, \beta, N)$ such that
$$|\tilde{T}_j^\beta f(x) - \tilde{T}_j^\beta f(z)| \leq C|x-z|^{N+\beta-\frac{N}{p}} \|f\|_{L^p(\mathbb{R}^N)}$$
for all $f \in L^p(\mathbb{R}^N)$ and $j = 1, \ldots, N$.

ii) If $p = \frac{N}{N-1+\beta}$ and $1 \leq q < p$, then there exists $C = C(q, \beta, N)$ such that
$$|\tilde{T}_j^\beta f(x) - \tilde{T}_j^\beta f(z)| \leq C|x-z| (|\ln|x-z||+1)^{\frac{1}{p'}} \left( \|f\|_{L^p(\mathbb{R}^N)} + \|f\|_{L^q(\mathbb{R}^N)} \right)$$
for all $f \in L^p(\mathbb{R}^N) \cap L^q(\mathbb{R}^N)$ and $j = 1, \ldots, N$.

**Remark 1.8** As we will see in Section 4, $\tilde{I}_N f = \tilde{T}^0 \cdot \mathcal{R} f$, where $\mathcal{R}$ is the vector Riesz transform. More generally, one could define modified Riesz potentials $\tilde{I}_{N+\beta} f$, for which one would have $\tilde{I}_{N+\beta} f = \tilde{T}^\beta \cdot \mathcal{R} f$. Therefore, the regularity properties of the maps $\tilde{T}$ imply regularity properties for the modified Riesz and logarithmic potentials.

Now, whereas our interest in the case $\beta = 0$ is related to the mapping properties of the logarithmic potential, and therefore regularity of solutions to Poisson's equation, we have further reason to study the mapping properties of $T_j^\beta$ for $\beta \in [0,1)$. Our motivation for considering this regime stems from the relationship of the Riesz potentials and Sobolev embeddings. In fact, one can show that Riesz potentials enable one to give a unified treatment of a number of Sobolev embeddings for $p > 1$. However, using Riesz potentials alone it is not possible to prove sharp embeddings concerning $H^{\alpha,p}(\mathbb{R}^N)$ for $\alpha \geq N$. In particular, this includes a spectrum of cases of the supercritical Sobolev embedding first proven by Brezis and Wainger [2][Corollary 5] in 1980. Among other results, they had shown the following theorem concerning the "almost" Lipschitz continuity of elements of the Bessel potential space $H^{\alpha,p}(\mathbb{R}^N)$.

**Remark 1.9** *Here, the notation $H^{\alpha,p}(\mathbb{R}^N)$ is used to denote the Bessel potential spaces, functions with the representation $u = g_\alpha * f$ for some $f \in L^p(\mathbb{R}^N)$, where $g_\alpha$ is the Bessel potential of order $\alpha$. In the case $\alpha = l \in \mathbb{N}$, it is a result of Calderón that $H^{\alpha,p}(\mathbb{R}^N) = W^{l,p}(\mathbb{R}^N)$ (see [20][Chapter 5, p.131 and p.165]).*

**Theorem 1.10** *Let $1 < p < +\infty$ and suppose $u \in H^{\alpha,p}(\mathbb{R}^N)$ with $\alpha = 1 + \frac{N}{p}$. Then there exists $C = C(p, N)$ such that*

$$|u(x) - u(z)| \leq C|x-z|\left(|\ln|x-z||+1\right)^{\frac{1}{p'}}\|u\|_{H^{\alpha,p}(\mathbb{R}^N)}.$$

In Section 4, we show how the potential estimates we have obtained can be used to deduce the preceding theorem in its full generality. This is accomplished via a representation for Sobolev functions via the Riesz potentials and Theorem A for the range $1 < \alpha < N$, while for the range $\alpha \in [N, N+1)$ we find an alternative representation via the Riesz-type potentials (1.6) and invoke Theorem B.

**Remark 1.11** *In the case $\alpha = N + 1$, $p = 1$, one has the stronger result that any $u \in W^{N+1,1}(\mathbb{R}^N)$ is Lipschitz, see Gagliardo's 1958 paper [5].*

The organization of the paper will be as follows. In Section 2, we will prove Theorems A and D, as well as a weaker version of Theorem B. In Section 3, we show the equivalence of $\tilde{I}_N f$ and $\tilde{T}^N \cdot \mathcal{R}f$, thereby deducing Theorems B and C. Finally, in Section 4 we show how our results can be used to deduce Theorem 1.10.

## 2 Potential Estimates

In this section, we prove potential estimates for the Riesz and Riesz-type potentials introduced in Section 1.2. Here, we restate the theorems for the convenience of the reader.

We begin with the following theorem, on the continuity of the Riesz potential in the supercritical case.

**Theorem A** *Let $1 < p < +\infty$ and $\alpha = 1 + \frac{N}{p} < N$. Then for any $1 \leq q < p$, there exists $C = C(q, \alpha, N)$ such that*

$$|\tilde{I}_\alpha f(x) - \tilde{I}_\alpha f(z)| \leq C|x-z|\left(|\ln|x-z||+1\right)^{\frac{1}{p'}}\left(\|f\|_{L^q(\mathbb{R}^N)} + \|f\|_{L^p(\mathbb{R}^N)}\right)$$

*for all $f \in L^p(\mathbb{R}^N) \cap L^q(\mathbb{R}^N)$.*

**Proof.** For brevity, we write $\gamma = \gamma(N, \alpha)$ and define $r := |x - z|$. We have

$$|\tilde{I}_\alpha f(x) - \tilde{I}_\alpha f(z)| \leq \gamma \int_{B(x,2r)} |f(y)| \left(\frac{1}{|x-y|^{N-\alpha}} + \frac{1}{|z-y|^{N-\alpha}}\right) dy$$

$$+ \gamma \int_{B(x,2r)^c} |f(y)| \left|\frac{1}{|x-y|^{N-\alpha}} - \frac{1}{|z-y|^{N-\alpha}}\right| dy$$

$$=: I + II.$$

For $I$, an application of Hölder's inequality, along with the fact that for all $\theta > 0$ and $z$ such that $r = |z - x|$, we have the inequality

$$\int_{B(x,2r)} \frac{1}{|z-y|^\theta}\, dy \leq \int_{B(x,2r)} \frac{1}{|x-y|^\theta}\, dy$$

allows us to deduce that

$$I \leq 2\gamma \|f\|_{L^p(\mathbb{R}^N)} \left( \int_{B(x,2r)} \frac{1}{|x-y|^{(N-\alpha)p'}}\, dy \right)^{\frac{1}{p'}}$$

$$= 2\gamma |S^{N-1}|^{\frac{1}{p'}} \|f\|_{L^p(\mathbb{R}^N)} \left( \int_0^{2r} t^{N-1-(N-\alpha)p'}\, dt \right)^{\frac{1}{p'}}$$

$$= \frac{4\gamma |S^{N-1}|^{\frac{1}{p'}}}{(p')^{\frac{1}{p'}}} r \|f\|_{L^p(\mathbb{R}^N)}$$

$$\leq \frac{4\gamma |S^{N-1}|^{\frac{1}{p'}}}{(p')^{\frac{1}{p'}}} r \left(|\ln r| + 1\right)^{\frac{1}{p'}} \|f\|_{L^p(\mathbb{R}^N)}.$$

In order to treat the term $II$, we define the auxiliary function

$$h(t) := \frac{1}{|tx + (1-t)z - y|^{N-\alpha}},$$

and observe that for $y \in B(x, 2r)^c$, $h : [0,1] \to \mathbb{R}$ is smooth. Thus, the mean value theorem implies

$$h(1) - h(0) = h'(t_0)$$

for some $t_0 = t_0(x, z) \in (0, 1)$. Rewriting this equality in terms of the function $h$ we have

$$\frac{1}{|x-y|^{N-\alpha}} - \frac{1}{|z-y|^{N-\alpha}}$$
$$= \frac{-N+\alpha}{|t_0 x + (1-t_0)z - y|^{N-\alpha+1}} \frac{(t_0 x + (1-t_0)z - y) \cdot (x-z)}{|t_0 x + (1-t_0)z - y|}.$$

From this we conclude that

$$\left| \frac{1}{|x-y|^{N-\alpha}} - \frac{1}{|z-y|^{N-\alpha}} \right| \leq |x-z| \frac{2^{N-\alpha+1}(N-\alpha)}{|x-y|^{N-\alpha+1}},$$

where we have used the fact that $|t_0 x + (1-t_0)z - y| \geq \frac{1}{2}|x-y|$ for any $t_0 \in (0,1)$.

Now, $f \in L^p(\mathbb{R}^N) \cap L^q(\mathbb{R}^N)$ implies $f \in L^s(\mathbb{R}^N)$ for all $s \in [q, p]$, and we can thus estimate

$$II \leq r \left(2^{N-\alpha+1}(N-\alpha)\gamma\right) \int_{B(x,2r)^c} |f(y)| \frac{1}{|x-y|^{N-\alpha+1}}\, dy$$

$$\leq r \left(2^{N-\alpha+1}(N-\alpha)\gamma\right) |S^{N-1}|^{\frac{1}{s'}} \|f\|_{L^s(\mathbb{R}^N)} \left( \int_{2r}^\infty \frac{1}{t^{(N-\alpha+1)s'}} t^{N-1}\, dt \right)^{\frac{1}{s'}}.$$

Then since Young's inequality implies

$$\|f\|_{L^s(\mathbb{R}^N)} \leq \|f\|_{L^q(\mathbb{R}^N)} + \|f\|_{L^p(\mathbb{R}^N)},$$

and $|S^{N-1}|^{\frac{1}{s'}} \leq \min\left\{|S^{N-1}|^{\frac{1}{p'}}, |S^{N-1}|^{\frac{1}{q'}}\right\}$, it only remains to show that we can in fact choose an $s = s(r) \in [q, p]$ such that

$$\left(\int_{2r}^\infty \frac{1}{t^{(N-\alpha+1)s'}} t^{N-1}\, dt\right)^{\frac{1}{s'}} \leq C(p,q,N)\left(|\ln r| + 1\right)^{\frac{1}{p'}},$$

and the result is demonstrated. First, we estimate the integral

$$\left(\int_{2r}^\infty \frac{1}{t^{(N-\alpha+1)s'}} t^{N-1}\, dt\right)^{\frac{1}{s'}} = \frac{(2r)^{N(\frac{1}{p} - \frac{1}{s})}}{(Ns'(\frac{1}{s} - \frac{1}{p}))^{\frac{1}{s'}}} \leq \frac{(2r)^{N(\frac{1}{p} - \frac{1}{s})}}{(N(\frac{1}{s} - \frac{1}{p}))^{\frac{1}{s'}}}.$$

Now, since $N(\frac{1}{p} - \frac{1}{q}) \in (-\infty, 0)$, there exists a unique $0 < \epsilon_0 = \epsilon_0(p,q) < 1$ such that $N(\frac{1}{p} - \frac{1}{q}) = \frac{-1}{|\ln \epsilon_0|}$. Then when $r \in [\epsilon_0, \infty)$, the choice of $s = q$ in the above integral allows us to obtain the bound

$$\frac{(2r)^{N(\frac{1}{p} - \frac{1}{q})}}{(N(\frac{1}{q} - \frac{1}{p}))^{\frac{1}{q'}}} = (2r)^{\frac{-1}{|\ln \epsilon_0|}} |\ln \epsilon_0|^{\frac{1}{q'}} \leq (2\epsilon_0)^{\frac{-1}{|\ln \epsilon_0|}} \left(|\ln \epsilon_0|\right)^{\frac{1}{q'}} \left(|\ln r| + 1\right)^{\frac{1}{p'}}.$$

Finally, we consider the case when $r \in (0, \epsilon_0)$. Here, we define $s = s(r)$ by $N(\frac{1}{p} - \frac{1}{s}) = \frac{-1}{|\ln r|}$. It follows from the definition of $\epsilon_0$ that $s \in (q, p)$. Utilizing this exponent in the computation for $II$ yields

$$\frac{(2r)^{N(\frac{1}{p} - \frac{1}{s})}}{(N(\frac{1}{s} - \frac{1}{p}))^{\frac{1}{s'}}} = 2^{\frac{-1}{|\ln r|}} r^{\frac{-1}{|\ln r|}} |\ln r|^{\frac{1}{s'}} \leq exp(1)\left(|\ln r| + 1\right)^{\frac{1}{p'}},$$

which implies the desired result. ∎

We now prove the following theorem concerning the mapping properties of the modified logarithmic potential. Although the estimates in this theorem are not as sharp as those claimed in Theorem B, we will see that the proof follows in analogy to Theorem A. Moreover, we will utilize this estimate in order to obtain the sharp result, which we prove in Section 3.

**Theorem B'** *Let $1 < p \leq \frac{N}{N-1}$.*

i) *If $1 < p < \frac{N}{N-1}$, then there exists $C = C(p, N)$ such that*

$$|\tilde{I}_N f(x) - \tilde{I}_N f(z)| \leq C|x - z|^{N - \frac{N}{p}} \left(|\ln |x - z|| + 1\right) \|f\|_{L^p(\mathbb{R}^N)}$$

*for all $f \in L^p(\mathbb{R}^N)$.*

ii) *If $p = \frac{N}{N-1}$ and $1 \leq q < p$, then there exists $C = C(q, N)$ such that*

$$|\tilde{I}_N f(x) - \tilde{I}_N f(z)| \leq C|x - z|\left(|\ln |x - z|| + 1\right)\left(\|f\|_{L^p(\mathbb{R}^N)} + \|f\|_{L^q(\mathbb{R}^N)}\right)$$

*for all $f \in L^p(\mathbb{R}^N) \cap L^q(\mathbb{R}^N)$.*

**Proof.** We begin under the hypothesis of part *i)*, assuming $f \in L^p(\mathbb{R}^N)$ for some $1 < p < \frac{N}{N-1}$. As in the proof of Theorem A, we split the estimate into two piece. We have for $r = |z - x|$,

$$|\tilde{I}_N f(x) - \tilde{I}_N f(z)| \leq \omega_N \int_{B(x,2r)} |f(y)|\, (|\ln|x-y|| + |\ln|z-y||)\; dy$$

$$+ \omega_N \int_{B(x,2r)^c} |f(y)|\, |\ln|x-y| - \ln|z-y||\; dy$$

$$=: I + II.$$

For $I$, we this time require a different estimate in several regimes of $r \in (0, +\infty)$. First notice that since $B(x, 2r) \subset B(z, 3r)$, we have

$$I \leq \omega_N \left( \int_{B(x,2r)} |f(y)|\, |\ln|x-y||\; dy + \int_{B(z,3r)} |f(y)|\, |\ln|z-y||\; dy \right)$$

$$\leq 2\omega_N |S^{N-1}|^{\frac{1}{p'}} \|f\|_{L^p(\mathbb{R}^N)} \left( \int_0^{3r} |\ln t|^{p'} t^{N-1} dt \right)^{\frac{1}{p'}}.$$

When $r \leq \epsilon_1$, for some $0 < \epsilon_1 < \frac{1}{3}$ to be chosen shortly, we proceed as follows. The change of variables $t = exp(-x)$ implies

$$\left( \int_0^{3r} |\ln t|^{p'} t^{N-1} dt \right)^{\frac{1}{p'}} = \left( \int_{|\ln 3r|}^{\infty} x^{p'} e^{-Nx}\; dx \right)^{\frac{1}{p'}},$$

for which we can use the following asymptotics of the upper incomplete gamma function

$$\frac{\int_w^\infty x^s e^{-Nx}\; dx}{w^s e^{-Nw}} \to 1$$

as $w \to \infty$ to deduce that there exists and $\epsilon_1 \in (0, \frac{1}{3})$ such that for $r \leq \epsilon_1$ we have

$$\left( \int_{|\ln 3r|}^\infty x^{p'} e^{-Nx}\; dx \right)^{\frac{1}{p'}} \leq 2|\ln 3r| e^{-\frac{N}{p'}|\ln 3r|}$$

$$= 2|\ln 3r|(3r)^{N-\frac{N}{p}}.$$

For $r \in [\epsilon_1, \frac{1}{3}]$, we have

$$\left( \int_0^{3r} |\log t|^{p'} t^{N-1} dt \right)^{\frac{1}{p'}} \leq \left( \int_0^1 |\ln t|^{p'} t^{N-1} dt \right)^{\frac{1}{p'}}$$

$$\leq C \left( \frac{r}{\epsilon_1} \right)^{N-\frac{N}{p}},$$

since $\frac{r}{\epsilon_1} \geq 1$. Finally, when $r \in (\frac{1}{3}, \infty)$, we have

$$\left(\int_0^{3r} |\log t|^{p'} t^{N-1} dt\right)^{\frac{1}{p'}} = \left(\int_0^1 |\ln t|^{p'} t^{N-1} dt + \int_1^{3r} |\ln t|^{p'} t^{N-1} dt\right)^{\frac{1}{p'}}$$

$$\leq \left(C + |\ln 3r|^{p'} \left(\frac{(3r)^N - 1}{N}\right)\right)^{\frac{1}{p'}}$$

$$\leq Cr^{N-\frac{N}{p}} (|\ln r| + 1).$$

Combining these estimates, we conclude that

$$I \leq C \|f\|_{L^p(\mathbb{R}^N)} r^{N-\frac{N}{p}} (|\ln r| + 1).$$

For $II$, we can proceed as before, since the mean value theorem again implies

$$|\log |x-y| - \log |z-y|| \leq 2|x-z|\frac{1}{|y-x|},$$

for all $y \in B(x, 2r)^c$ and therefore,

$$II \leq 2r|S^{N-1}|^{\frac{1}{p'}} \|f\|_{L^p(\mathbb{R}^N)} \left(\int_{2r}^{\infty} \frac{1}{t^{p'}} t^{N-1} dt\right)^{\frac{1}{p'}}.$$

Part *i)* then follows, since

$$\left(\int_{2r}^{\infty} \frac{1}{t^{p'}} t^{N-1} dt\right)^{\frac{1}{p'}} = \frac{(2r)^{\frac{N}{p'}-1}}{(p'-N)^{\frac{1}{p'}}}.$$

For the proof of Part *ii)*, the preceding analysis for $I$ is exactly the same while $II$ can be estimated as that in the proof of Theorem A. ∎

Finally, we conclude this section with a proof of Theorem D, which follows in analogy to the proofs of the previous two theorems.

**Theorem D** *Let $\beta \in [0, 1)$.*

*i) If $1 < p < \frac{N}{N-1+\beta}$, then there exists $C = C(p, \beta, N)$ such that*

$$|\tilde{T}_j^\beta f(x) - \tilde{T}_j^\beta f(z)| \leq C|x-z|^{N+\beta-\frac{N}{p}} \|f\|_{L^p(\mathbb{R}^N)}$$

*for all $f \in L^p(\mathbb{R}^N)$ and $j = 1, \ldots, N$.*

*ii) If $p = \frac{N}{N-1+\beta}$ and $1 \leq q < p$, then there exists $C = C(q, \beta, N)$ such that*

$$|\tilde{T}_j^\beta f(x) - \tilde{T}_j^\beta f(z)| \leq C|x-z| (|\ln |x-z|| + 1)^{\frac{1}{p'}} (\|f\|_{L^p(\mathbb{R}^N)} + \|f\|_{L^q(\mathbb{R}^N)})$$

*for all $f \in L^p(\mathbb{R}^N) \cap L^q(\mathbb{R}^N)$ and $j = 1, \ldots, N$.*

**Proof.** Part *i)* can be argued in a similar manner as in the proofs of Theorems A and B', this time applying Hölder's inequality for $f \in L^p(\mathbb{R}^N)$ with $1 < p < \frac{N}{N-1}$ for both $I$ and $II$.

The proof of Part *ii)* is analogous to that of Theorem A. While $I$ is estimated as in Part *i)*, for $II$ we observe that the mean value theorem implies that for any $y \in B(x, 2r)^c$

$$\left|\frac{x_j - y_j}{|x-y|^{1-\beta}} - \frac{z_j - y_j}{|z-y|^{1-\beta}}\right| \leq (2-\beta)2^{1-\beta}|x-z|\frac{1}{|x-y|^{1-\beta}},$$

and the log-Lipschitz estimate follows directly from the remainder of the argument of Theorem A. ∎

## 3 Regularity for Solutions to Poisson's Equation

Before we prove Theorem B, we first record the following lemma connecting the potentials $\tilde{T}^N$ and the Riesz potentials (1.1). Let us remark that this generalizes the case $N = 2$, which has been proven in [6].

**Lemma 3.1** *Suppose $f \in C_c^\infty(\mathbb{R}^N)$. Then*
$$-\Delta\left(\tilde{T}^0 \cdot \mathcal{R}f\right) = I_{N-2}f,$$
*with the convention that $I_0 f = f$.*

**Proof.** It suffices to show that
$$-\int_{\mathbb{R}^N} (\tilde{T}^0 \cdot \mathcal{R}f)\Delta\varphi = -\int_{\mathbb{R}^N} (I_{N-2}f)\varphi$$
for all $\varphi \in C_c^\infty(\mathbb{R}^N)$. First, we observe that in the proof of Theorem D, one has the stronger inequality
$$\int_{\mathbb{R}^N} \left|\frac{y_j - x_j}{|y - x|} - \frac{y_j}{|y|}\right| |R_j f(y)|\, dy \leq C|x|^{N-\frac{N}{p}} \|R_j f\|_{L^p(\mathbb{R}^N)}$$
$$\leq C|x|^{N-\frac{N}{p}} \|f\|_{L^p(\mathbb{R}^N)}$$
for $j = 1, 2, \ldots, N$ and $1 < p < \frac{N}{N-1}$. Therefore, we may utilize Fubini's theorem to deduce
$$-\int_{\mathbb{R}^N} (\tilde{T}^0 \cdot \mathcal{R}f)\Delta\varphi$$
$$= \frac{-\gamma(N, N-1)}{N-1} \sum_{j=1}^N \int_{\mathbb{R}^N} \left(\int_{\mathbb{R}^N} \left[\frac{y_j - x_j}{|y - x|} - \frac{y_j}{|y|}\right] R_j f(y)\, dy\right) \Delta\varphi(x)\, dx$$
$$= \frac{-\gamma(N, N-1)}{N-1} \sum_{j=1}^N \int_{\mathbb{R}^N} \left(\int_{\mathbb{R}^N} \left[\frac{y_j - x_j}{|y - x|} - \frac{y_j}{|y|}\right] \Delta\varphi(x)\, dx\right) R_j f(y)\, dy.$$

Further, since the divergence theorem implies $\int_{\mathbb{R}^N} \Delta\varphi(x)\, dx = 0$, we have that
$$\int_{\mathbb{R}^N} \left[\frac{y_j - x_j}{|y - x|} - \frac{y_j}{|y|}\right] \Delta\varphi(x)\, dx = \int_{\mathbb{R}^N} \frac{y_j - x_j}{|y - x|} \Delta\varphi(x)\, dx.$$

Now, we define
$$g_j(y) := \frac{-\gamma(N, N-1)}{N-1} \int_{\mathbb{R}^N} \frac{y_j - x_j}{|y - x|} \Delta\varphi(x)\, dx.$$

If we can show that $g_j = R_j(I_{N-2}\varphi)$ almost everywhere, then we would have
$$-\int_{\mathbb{R}^N} \tilde{T}^0 \cdot \mathcal{R}f \Delta\varphi = \sum_{j=1}^N \int_{\mathbb{R}^N} R_j(I_{N-2}\varphi) R_j f$$
$$= \int_{\mathbb{R}^N} f(I_{N-2}\varphi)$$
$$= \int_{\mathbb{R}^N} (I_{N-2}f)\varphi,$$

which is the thesis. Notice that

$$g_j(y) = \frac{-y_j}{N-1}(I_{N-1}\Delta\varphi)(y) + \frac{1}{N-1}I_{N-1}(x_j\Delta\varphi)(y),$$

and therefore,

$$\begin{aligned}\widehat{g_j}(\xi) &= \frac{1}{2\pi i}\frac{1}{N-1}\frac{\partial}{\partial\xi_j}\left((2\pi|\xi|)^{1-N}\widehat{\Delta\varphi}(\xi)\right) + \frac{1}{N-1}(2\pi|\xi|)^{1-N}\left(\widehat{x_j\Delta\varphi}(\xi)\right)\\ &= \frac{i(2\pi)^{2-N}}{N-1}\left[\frac{\partial}{\partial\xi_j}(|\xi|^{3-N}\widehat{\varphi}(\xi)) - \frac{1}{|\xi|^{N-1}}\frac{\partial}{\partial\xi_j}(|\xi|^2\widehat{\varphi}(\xi))\right]\\ &= -i\frac{\xi_j}{|\xi|}(2\pi|\xi|)^{2-N}\widehat{\varphi}(\xi)\\ &= (R_j I_{N-2}\varphi)\widehat{\phantom{.}}.\end{aligned}$$

Here, the above should be interpreted in the sense of tempered distributions, and with the convention

$$\widehat{\varphi}(\xi) = \int_{\mathbb{R}^N}\varphi(x)e^{-2\pi i x\cdot\xi}\,dx$$

for the Fourier transform. Thus, we have proved that $g_j = R_j I_{N-2}\varphi$ as distributions, which implies almost everywhere equality as functions, and the result is demonstrated. ∎

Again, for the convenience of the reader we here recall the statement of Theorem B before giving its proof.

**Theorem B** *Let* $1 < p \leq \frac{N}{N-1}$.

  i) *If* $1 < p < \frac{N}{N-1}$, *then there exists* $C = C(p, N)$ *such that*

$$|\tilde{I}_N f(x) - \tilde{I}_N f(z)| \leq C|x-z|^{N-\frac{N}{p}}\|f\|_{L^p(\mathbb{R}^N)}$$

  *for all* $f \in L^p(\mathbb{R}^N)$.

  ii) *If* $p = \frac{N}{N-1}$ *and* $1 \leq q < p$, *then there exists* $C = C(q, N)$ *such that*

$$|\tilde{I}_N f(x) - \tilde{I}_N f(z)| \leq C|x-z|(|\ln|x-z||+1)^{\frac{1}{N}}(\|f\|_{L^p(\mathbb{R}^N)} + \|f\|_{L^q(\mathbb{R}^N)})$$

  *for all* $f \in L^p(\mathbb{R}^N) \cap L^q(\mathbb{R}^N)$.

**Proof.** We will show that

$$\tilde{I}_N f \equiv \tilde{T}^0 \cdot \mathcal{R}f \tag{3.1}$$

for all $f \in L^p(\mathbb{R}^N)$, $1 < p < \frac{N}{N-1}$. The claimed regularity estimate for $\tilde{I}_N f$ then follows from Theorem D and boundedness of $\mathcal{R} : L^p(\mathbb{R}^N) \to L^p(\mathbb{R}^N)$ for $1 < p < +\infty$.

First, we observe that it suffices to establish (3.1) for $f \in C_c^\infty(\mathbb{R}^N)$, since an application of Theorems B′ and D yields the estimates

$$|\tilde{I}_N(f_n - f)(x)| \leq C|x|^{N-\frac{N}{p}}(|\ln|x||+1)\|f_n - f\|_{L^p(\mathbb{R}^N)}$$
$$|\tilde{T}^0 \cdot \mathcal{R}(f_n - f)(x)| \leq C|x|^{N-\frac{N}{p}}\|f_n - f\|_{L^p(\mathbb{R}^N)}.$$

One is then able to conclude (3.1) for all $f \in L^p(\mathbb{R}^N)$, using the fact that $\tilde{I}_N f$, $\tilde{T}^0 \cdot \mathcal{R}f$ are continuous for $f \in L^p(\mathbb{R}^N)$, $1 < p < \frac{N}{N-1}$.

We therefore establish (3.1). Now, we already know that if $f \in C_c^\infty(\mathbb{R}^N)$ that

$$-\Delta \tilde{I}_N f = I_{N-2} f.$$

Further, Lemma 3.1 asserts that

$$-\Delta \tilde{T}^0 \cdot \mathcal{R}f = I_{N-2}f,$$

and so if we define $w := \tilde{I}_N f - \tilde{T}^0 \cdot \mathcal{R}f$, then we have that $w$ is harmonic. Moreover, Theorems B$'$ and D imply that

$$\frac{w(x)}{|x|} \to 0 \text{ as } |x| \to \infty.$$

But then one can invoke Liouville's theorem to conclude that $w \equiv c$, and since $w(0) = 0$ this implies $w \equiv 0$. This proves the claim and hence the theorem. ∎

As a consequence of Theorem 1.3 and Theorems A, B and D, we deduce the uniform regularity estimates for $\tilde{I}_2$, and therefore obtain Theorem C.

## 4 Sobolev Embedding Theorem in the Supercritical Exponent

In this section, we show how the preceding results for potentials, combined with appropriate representations for Sobolev functions, can be used to deduce Theorem 1.10.

The idea for the proof is as follows. Observe that for $u \in C_c^\infty(\mathbb{R}^N)$ and $0 < m < N$ even, one has the potential representation

$$u = I_m (-\Delta)^{\frac{m}{2}} u.$$

Thus, if $m = 1 + \frac{N}{p}$, one can apply Theorem A to deduce the inequality

$$|u(x) - u(z)| \le C|x - z|(|\ln|x-z|| + 1)^{\frac{1}{p'}} \left( \left\|(-\Delta)^{\frac{m}{2}} u\right\|_{L^q(\mathbb{R}^N)} + \left\|(-\Delta)^{\frac{m}{2}} u\right\|_{L^p(\mathbb{R}^N)} \right).$$

Then applying Holder's inequality, and using the fact that $u$ has compact support, one obtains the inequality

$$|u(x) - u(z)| \le C|x - z|(|\ln|x-z|| + 1)^{\frac{1}{p'}} \left\|(-\Delta)^{\frac{m}{2}} u\right\|_{L^p(\mathbb{R}^N)}.$$

Now, here $C$ depends on the support of $u$, and so one cannot use density to conclude the preceding inequality holds for general $u \in W^{m,p}(\mathbb{R}^N)$. Moreover, the handling of the case $m$ even or odd, let alone non-integer as in the setting of Bessel potential spaces make the matter slightly more complicated. This clearly tells us that a more delicate analysis is required to obtain the result, which we will now prove in its full generality.

**Proof of Theorem 1.10.** We claim that it suffices to show that for any $u \in C_c^\infty(\mathbb{R}^N)$ we have the estimate

$$|u(x) - u(z)| \leq C|x - z|(|\ln|x - z|| + 1)^{\frac{1}{p'}} \|u\|_{H^{\alpha,p}(\mathbb{R}^N)} \quad (4.1)$$

for all $x, z \in \mathbb{R}^N$ such that $|x - z| \leq 1$. If so, then this implies the inequality

$$|u(x)| \leq |u(z)| + |u(x) - u(z)|$$
$$\leq |u(z)| + C|x - z|(|\ln|x - z|| + 1)^{\frac{1}{p'}} \|u\|_{H^{\alpha,p}(\mathbb{R}^N)}.$$

Then averaging over $z \in B(x, 1)$, and using Hölder's inequality we have

$$|u(x)| \leq C\|u\|_{H^{\alpha,p}(\mathbb{R}^N)},$$

so that taking the supremum over all $x \in \mathbb{R}^N$, we deduce that

$$\|u\|_{L^\infty(\mathbb{R}^N)} \leq C\|u\|_{H^{\alpha,p}(\mathbb{R}^N)}$$

for all $u \in C_c^\infty(\mathbb{R}^N)$. Thus, if $|x - z| \geq 1$, we have

$$|u(x) - u(z)| \leq 2|x - z|(|\ln|x - z|| + 1)^{\frac{1}{p'}} \|u\|_{L^\infty(\mathbb{R}^N)},$$

and so if (4.1) holds for $|x-z| \leq 1$, it continues to hold for all $x, z \in \mathbb{R}^N$, perhaps increasing the constant. The result then follows for general $u \in H^{\alpha,p}(\mathbb{R}^N)$ by density, since the above inequality shows that any sequence Cauchy in $H^{\alpha,p}(\mathbb{R}^N)$ is Cauchy in the space $Lip^{0,-\frac{1}{p'}}(\mathbb{R}^N)$, defined by

$$\|u\|_{Lip^{0,-\frac{1}{p'}}(\mathbb{R}^N)} := \|u\|_{L^\infty(\mathbb{R}^N)} + \sup_{x \neq z} \frac{|u(x) - u(z)|}{|x - z|(|\ln|x - z|| + 1)^{\frac{1}{p'}}}.$$

**Case 1:** $1 < \alpha < N$.

To establish (4.1), let $\alpha = m + 2\gamma$ for $m \in \mathbb{N}$ and $\gamma \in [0, \frac{1}{2})$. Let $\phi \in C_c^\infty(B(0, 4))$ be a cutoff function such that

$$\phi \equiv 1 \text{ on } B(0, 2);$$
$$\phi \equiv 0 \text{ outside } B(0, 4);$$
$$\|D^\nu \phi\|_{L^\infty(B(0,4))} \leq C;$$
$$\|(-\Delta)^\gamma D^\nu \phi\|_{L^\infty(\mathbb{R}^N)} \leq C,$$

for all $0 \leq |\nu| \leq m$. Given $u \in C_c^\infty(\mathbb{R}^N)$ and $x \in \mathbb{R}^N$, define $h := u\phi_x = u(\cdot)\phi(\cdot - x)$. Now, as $h \in C_c^\infty(B(x, 4)) \subset C_c^\infty(\mathbb{R}^N)$, one can verify by taking the Fourier transform that $h$ has the representation $h = I_\alpha g$ for

$$g := \begin{cases} (-\Delta)^\gamma (-\Delta)^{\frac{m}{2}} h & \text{; if } m \text{ is even} \\ \sum_{i=1}^N R_j(-\Delta)^\gamma \frac{\partial}{\partial x_j}(-\Delta)^{\frac{m-1}{2}} h & \text{; if } m \text{ is odd.} \end{cases} \quad (4.2)$$

Now, since $h \in C_c^\infty(\mathbb{R}^N)$, it is easy to verify that $g \in L^q(\mathbb{R}^N)$ for all $1 < q \leq \infty$. Thus, an application of Theorem A to $g$ for any $1 < q < p$ implies the inequality

$$|I_\alpha g(x) - I_\alpha g(y)| \leq C|x - z|(|\ln|x - z|| + 1)^{\frac{1}{p'}} (\|g\|_{L^q(\mathbb{R}^N)} + \|g\|_{L^p(\mathbb{R}^N)}),$$

and so it suffices to prove that there exists $C$ depending only on $\phi$ such that

$$\|g\|_{L^q(\mathbb{R}^N)} + \|g\|_{L^p(\mathbb{R}^N)} \leq C\|u\|_{H^{\alpha,p}(\mathbb{R}^N)},$$

and we will have demonstrated the inequality (4.1). If $\gamma = 0$, this is a consequence of the boundedness of the Riesz transform on $L^q(\mathbb{R}^N)$ and $L^p(\mathbb{R}^N)$, Hölder's inequality to bound the $L^q(B(x,4))$ norm by the $L^p(B(x,4))$ norm, and the product rule for derivatives of $u\phi_x$.

When $\gamma \in (0, \frac{1}{2})$, the result still holds through slightly more subtle analysis. After using boundedness of the Riesz transform, we here use the analogy of Hölder's inequality for the fractional Laplacian

$$\|(-\Delta)^\gamma w\|_{L^q(\mathbb{R}^N)} \leq C \|(-\Delta)^\gamma w\|_{L^p(\mathbb{R}^N)},$$

for $w = (-\Delta)^{\frac{m}{2}} h$ or $w = \frac{\partial}{\partial x_j}(-\Delta)^{\frac{m-1}{2}} h$ depending on parity, which can be found in [18][Lemma A.1, p. 46], where one can check that $C = C(diam(supp\, h))$, and therefore can be taken depending only on $\phi$. The last result to check is whether there is an analogy for the product rule, from which we could conclude

$$\|(-\Delta)^\gamma w\|_{L^p(\mathbb{R}^N)} \leq C\|u\|_{H^{\alpha,p}(\mathbb{R}^N)}.$$

If we observe that

$$w = \sum_{|\mu|+|\nu|=m} C(\mu,\nu) D^\mu u D^\nu \phi_x,$$

then it suffices to show that

$$\|(-\Delta)^\gamma (D^\mu u D^\nu \phi_x)\|_{L^p(\mathbb{R}^N)} \leq C\|u\|_{H^{\alpha,p}(\mathbb{R}^N)}$$

for any $|\mu|+|\nu| = m$. Now, using the singular integral definition of the fractional Laplacian (which is not, in fact, singular, since $u, \phi \in C_c^\infty(\mathbb{R}^N)$ and $\gamma \in (0, \frac{1}{2})$), it suffices to obtain bounds for

$$A := \|(-\Delta)^\gamma D^\mu u D^\nu \phi_x\|_{L^p(\mathbb{R}^N)}.$$

However, we can estimate

$$A^p = c_{N,\gamma} \int_{\mathbb{R}^N} \left| \int_{\mathbb{R}^N} \frac{(D^\mu u D^\nu \phi_x)(z) - (D^\mu u D^\nu \phi_x)(y)}{|z-y|^{N+2\gamma}} \, dy \right|^p dz$$

$$\leq 2^{p-1} c_{N,\gamma} \int_{\mathbb{R}^N} |D^\nu \phi(z-x)|^p \left| \int_{\mathbb{R}^N} \frac{D^\mu u(z) - D^\mu u(y)}{|z-y|^{N+2\gamma}} \, dy \right|^p dz$$

$$+ 2^{p-1} c_{N,\gamma} \int_{\mathbb{R}^N} \left| \int_{\mathbb{R}^N} D^\mu u(y) \frac{D^\nu \phi(z-x) - D^\nu \phi(y-x)}{|z-y|^{N+2\gamma}} \, dy \right|^p dz$$

$$=: I + II.$$

To estimate $I$, we use the fact that $\|D^\nu \phi\|_{L^\infty(\mathbb{R}^N)} \leq C$ to obtain the upper bound

$$I \leq C \|(-\Delta)^\gamma D^\mu u\|_{L^p(\mathbb{R}^N)}^p,$$

while we break the integral in $II$ into four pieces and estimate them separately. First, we split the outer integral into two pieces,

$$II = \int_{B(x,6)} \left| \int_{\mathbb{R}^N} D^\mu u(y) \frac{D^\nu \phi(z-x) - D^\nu \phi(y-x)}{|z-y|^{N+2\gamma}} \, dy \right|^p dz$$

$$+ \int_{B(x,6)^c} \left| \int_{\mathbb{R}^N} D^\mu u(y) \frac{D^\nu \phi(z-x) - D^\nu \phi(y-x)}{|z-y|^{N+2\gamma}} \, dy \right|^p dz$$

$$=: II_1 + II_2.$$

We further break the inner integral in $II_1$ into two pieces, which yields the upper bound

$$II_1 \leq 2^{p-1} \int_{B(x,6)} \left| \int_{B(x,8)} D^\mu u(y) \frac{D^\nu \phi(z-x) - D^\nu \phi(y-x)}{|z-y|^{N+2\gamma}} \, dy \right|^p dz$$

$$+ 2^{p-1} \int_{B(x,6)} \left| \int_{B(x,8)^c} D^\mu u(y) \frac{D^\nu \phi(z-x) - D^\nu \phi(y-x)}{|z-y|^{N+2\gamma}} \, dy \right|^p dz.$$

Now, for the first term, we use the fact that $\phi$ is Lipschitz to bound

$$\int_{B(x,6)} \left| \int_{B(x,8)} D^\mu u(y) \frac{D^\nu \phi(z-x) - D^\nu \phi(y-x)}{|z-y|^{N+2\gamma}} \, dy \right|^p dz$$

$$\leq \int_{B(x,8)} \left( \int_{B(x,8)} \frac{|D^\mu u(y)|}{|z-y|^{N+2\gamma-1}} \, dy \right)^p dz.$$

Then using the fact that $\gamma \in (0, \frac{1}{2})$, we have $2\gamma - 1 < 0$, and so we can use the mapping properties of the restricted Riesz potential

$$I^\Omega_{1-2\gamma} f(z) := c \int_\Omega \frac{f(y)}{|z-y|^{N+2\gamma-1}} \, dy$$

for $\Omega = B(x,8)$ (c.f. [7], Chapter 7, Lemma 7.12, p. 152) to obtain that

$$\int_{B(x,8)} \left( \int_{B(x,8)} \frac{|D^\mu u(y)|}{|z-y|^{N+2\gamma-1}} \, dy \right)^p dz \leq C \|D^\mu u\|^p_{L^p(\mathbb{R}^N)},$$

and this completes the estimate for this piece. Returning to the second term of the upper bound for $II_1$, we have by Hölder's inequality

$$\int_{B(x,6)} \left| \int_{B(x,8)^c} D^\mu u(y) \frac{D^\nu \phi(z-x) - D^\nu \phi(y-x)}{|z-y|^{N+2\gamma}} \, dy \right|^p dz$$

$$= \int_{B(x,6)} \left| \int_{B(x,8)^c} D^\mu u(y) \frac{D^\nu \phi(z-x)}{|z-y|^{N+2\gamma}} \, dy \right|^p dz$$

$$\leq C \|D^\mu u\|^p_{L^p(\mathbb{R}^N)},$$

and thus the estimate for $II_1$ is complete. We therefore return to estimate $II_2$,

and break the inner integral into two pieces in a similar manner,

$$II_2 \leq \int_{B(x,6)^c} \left| \int_{B(x,4)} D^\mu u(y) \frac{D^\nu \phi(z-x) - D^\nu \phi(y-x)}{|z-y|^{N+2\gamma}} \, dy \right|^p dz$$

$$+ \int_{B(x,6)^c} \left| \int_{B(x,4)^c} D^\mu u(y) \frac{D^\nu \phi(z-x) - D^\nu \phi(y-x)}{|z-y|^{N+2\gamma}} \, dy \right|^p dz.$$

Here, we are able to use Minkowski's inequality for integrals followed by Hölder's inequality to estimate the first term, while the second term is identically zero, since $supp \, \phi \subset B(0,4)$. We therefore conclude that

$$A^p \leq C \left( \|(-\Delta)^\gamma D^\mu u\|_{L^p(\mathbb{R}^N)}^p + \|D^\mu u\|_{L^p(\mathbb{R}^N)}^p \right),$$

which is the claim and completes the proof in this regime.

**Case 2:** $\alpha = N + \beta; \beta \in [0, 1)$.

In this regime of $\alpha$, we will use an alternative representation for $h$, in terms of $g$ and the unmodified potentials $T_j^\beta$,

$$T_j^\beta f(x) := \frac{\gamma(N, N-1+\beta)}{N-1+\beta} \int_{\mathbb{R}^N} \frac{y_j - x_j}{|y-x|^{1-\beta}} f(y) \, dy.$$

We begin with the equality $h = I_{N-1+\beta} g$, for $g$ as in (4.2), depending on the parity of $m$. Then we claim $g \in L^\infty(\mathbb{R}^N)$ and $|g(y)| \leq \frac{C}{|y|^{N+\beta}}$. If $m$ is even, these bounds are a consequence of the principle value representation for the fractional Laplacian of a function with compact support. When $m$ is odd, we observe that for $\gamma \in [0, \frac{1}{2})$, one can perform an integration by parts to deduce the alternative representation

$$R_j(-\Delta)^\gamma \phi(x) = I_{1-2\gamma} \frac{\partial \phi}{\partial x_j}(x)$$

$$= c \int_{\mathbb{R}^N} \frac{\phi(x) - \phi(y)}{|x-y|^{N+2\gamma}} \cdot \frac{x_j - y_j}{|x-y|} \, dy.$$

Thus, the composition of a fractional Laplacian and a Riesz transform on function with compact support has $L^\infty$ and decay bounds whose character is in line with the fractional Laplacian of a function with compact support. Now, if we define $A_\epsilon := B(x, \frac{1}{\epsilon}) \setminus B(x, \epsilon)$, then we have

$$I_{N-1+\beta} g(x) = \gamma(N, N-1+\beta) \int_{\mathbb{R}^N} \frac{1}{|x-y|^{1-\beta}} g(y) \, dy$$

$$= \gamma(N, N-1+\beta) \lim_{\epsilon \to 0} \int_{A_\epsilon} \frac{1}{|x-y|^{1-\beta}} g(y) \, dy$$

$$= \frac{-\gamma(N, N-1+\beta)}{N-1+\beta} \lim_{\epsilon \to 0} \int_{A_\epsilon} div_y \left( \frac{x-y}{|x-y|^{1-\beta}} \right) g(y) \, dy.$$

Then we integrate by parts to obtain

$$\int_{A_\epsilon} div_y \left( \frac{x-y}{|x-y|^{1-\beta}} \right) g(y)\, dy = \int_{\partial B(x,\epsilon)} \frac{x-y}{|x-y|^{1-\beta}} \cdot n(y)g(y)\, dy$$
$$+ \int_{\partial B(x,\frac{1}{\epsilon})} \frac{x-y}{|x-y|^{1-\beta}} \cdot n(y)g(y)\, dy$$
$$- \int_{A_\epsilon} \frac{x-y}{|x-y|^{1-\beta}} \cdot \nabla g(y)\, dy.$$

The fact that $g \in L^\infty(\mathbb{R}^N)$ and the decay $|g(y)| \leq \frac{C}{|y|^{N+\beta}}$ (as discussed earlier) imply that the two boundary integrals vanish as $\epsilon \to 0$. We therefore conclude that

$$I_{N-1+\beta}g(x) = \frac{-\gamma(N, N-1+\beta)}{N-1+\beta} \lim_{\epsilon \to 0} \int_{A_\epsilon} \frac{y-x}{|x-y|^{1-\beta}} \cdot \nabla g(y)\, dy. \quad (4.3)$$

This shows that the limit on the right hand side of (4.3) exists for all $x \in \mathbb{R}^N$. As such, if we let $x, z \in \mathbb{R}^N$, then we have that

$$h(x) - h(z) = \frac{-\gamma(N, N-1+\beta)}{N-1+\beta} \lim_{\epsilon \to 0} \int_{A_\epsilon} \left[ \frac{y-x}{|x-y|^{1-\beta}} - \frac{y-z}{|z-y|^{1-\beta}} \right] \cdot \nabla g(y)\, dy$$
$$= \frac{-\gamma(N, N-1+\beta)}{N-1+\beta} \int_{\mathbb{R}^N} \left[ \frac{y-x}{|x-y|^{1-\beta}} - \frac{y-z}{|z-y|^{1-\beta}} \right] \cdot \nabla g(y)\, dy.$$

Here we have used Lebesgue's dominated convergence theorem to pass the limit in the last integral, since the integrand is absolutely integrable over $\mathbb{R}^N$ (which is a consequence of the estimates in the proof of Theorem D). This then says that

$$|h(x) - h(z)| \leq \sum_{j=1}^{N} \left| \tilde{T}_j^\beta \frac{\partial g}{\partial x_j}(x) - \tilde{T}_j^\beta \frac{\partial g}{\partial x_j}(z) \right|,$$

and so we may now argue analagously to Case 1, this time applying Theorem D to obtain the desired estimate. ∎

## Acknowledgements

The first author is supported in part by a research grant (No: 471/13) of Amos Nevo from the Israel Science Foundation and a postdoctoral fellowship from the Planning and Budgeting Committee of the Council for Higher Education of Israel. The second author is supported in part by a Technion Fellowship. The authors would like to thank Yehuda Pinchover, Georgios Psaradakis, Simeon Reich, Sundaram Thangavelu and Igor Verbitsky for their helpful comments during the preparation of this work.